\documentclass[11pt, english, draft]{article}
\usepackage{amssymb,amsmath,amsthm,tikz}
\usepackage{graphicx,color}
\usetikzlibrary{shapes.geometric}
\usetikzlibrary{graphs}
\usetikzlibrary{graphs.standard}
\newtheorem{thm}{Theorem}[section]

\newtheorem{cor}[thm]{Corollary}

\newtheorem{prop}[thm]{Proposition}
\newtheorem{remark}[thm]{Remark}
\newtheorem{lem}[thm]{Lemma}

\date{}

\begin{document}
\title{\bf Some regular signed graphs with only two distinct eigenvalues}
\author{\large   Farzaneh Ramezani
\\ {\it Department of Mathematics,}
\\ {\it  K.N.Toosi University of Technology, Tehran, Iran}
\\ {\it P.O. Box 16315-1618}}
\maketitle  \footnotetext[1]{\tt Email: ramezani@kntu.ac.ir.}

\noindent {\bf Abstract}.
We consider signed graphs, i.e, graphs with positive or negative
signs on their edges. We determine the admissible parameters for the $\{5,6,\ldots,10\}$-regular signed graphs which have only two distinct eigenvalues. For each obtained parameter we provide some examples of signed graphs having two distinct eigenvalues. It turns out to construction of infinitely many signed graphs of each mentioned valency with only two distinct eigenvalues. We prove that for any $k\geq 5$ there are infinitely many connected signed $k$-regular graphs having maximum eigenvalue $\sqrt{k}$. Moreover for each $m\geq 4$ we construct a signed $8$-regular graph with spectrum $[4^m,-2^{2m}]$. These yield infinite family of $k$-regular Ramanujan graphs, for each $k$.



\section{Introduction}
We consider only simple graphs, i.e., graphs with out loops and
multiple edges. A \textit{signature}
on a graph $\Gamma=(V_{\Gamma},E_{\Gamma})$ is a function $\sigma:E_{\Gamma}\rightarrow \{1,-1\}$. A graph
$\Gamma$ provided with a signature $\sigma$ is called a \textit{signed
graph}, and will be denoted by $\Sigma=(\Gamma,\sigma)$. We call the graph $\Gamma$ the
\textit{ground graph} of the signed graph $\Sigma$, and denote it by $|\Sigma|$. Two signed graphs $\Sigma_1$ and $\Sigma_2$
are \textit{isomorphic}, denoted by $\Sigma_1\cong \Sigma_2$, if there is a graph isomorphism from $|\Sigma_1|$ to $|\Sigma_2|$, which preserves the signs of edges. For a subset $X$ of vertices of a signed graph $\Sigma=(\Gamma,\sigma)$, the signed graph of $\Sigma$ \textit{induced} on $X$ is a signed graph where the ground is the subgraph of $|\Sigma|$ induced on $X$ where the corresponding sign function is the restriction of $\sigma$ to $E_{\langle X\rangle}$. By $\Sigma\setminus X$ we mean the signed graph induced on $V_{|\sum|}\setminus X$. With a \textit{signed $k$-regular graph}, we mean a signed graph whose ground is $k$-regular. For a graph $\Gamma$, by $\Gamma^+$ (resp. $\Gamma^-$) we mean the all positive (resp. the all negative) signed graph with ground $\Gamma$. For a signed graph $\Sigma=(\Gamma,\sigma)$, by $-\Sigma$, we mean the signed graph $(\Gamma,-\sigma)$.

The \textit{adjacency matrix}, $A_{\Sigma}$ of the signed graph $\Sigma=(\Gamma,\sigma)$ on
the vertex set $V=\{v_1,v_2,\ldots,v_n\}$, is an $n\times n$
matrix whose entries are $$A_{\Sigma}(i,j)=\left\{
\begin{array}{ll}
\sigma(v_iv_j), & \hbox{if $v_i$ is adjacent to $v_j$;} \\
0, & \hbox{otherwise.}
\end{array}
\right.
$$ The \textit{ordinary adjacency matrix} of the graph $\Gamma$, is denoted by $A_{\Gamma}$. Whose entries
are the absolute values of the entries of $A_{\Sigma}$.

The \textit{spectrum} of a signed graph is the eigenvalues of its
adjacency matrix and will be denoted by $[\lambda_1^{m_1},\lambda_2^{m_2},\ldots,\lambda_s^{m_s}]$, where $m_i$ is the multiplicity of eigenvalue $\lambda_i$, for $i=1,2,\ldots,s$. By $m(\lambda)=m_{\Sigma}(\lambda)$, we denote the multiplicity of $\lambda$ as an eigenvalue of the signed graph $\Sigma$. By $O_n$, $J_n$ and $I_n$, we mean the all zero
matrix, all one matrix, and identity matrix of order $n$,
respectively. If $n$ is already known, then we
simply use $O,J,I$. For a matrix $X$ of order $m\times n$, by $X_i$ we mean the $i$th row of $X$, for $i=1,2,\ldots,m$.

An $n\times n$ matrix $C$, with $(0,\pm 1)$-entries is called a \textit{weighing matrix}, if $CC^t=C^tC=\alpha I_n$, where $\alpha$ is a positive integer, called the \textit{weight} of the weighing matrix $C$. A signed graph is called \textit{weighing} if its adjacency matrix is a weighing matrix. Weighing matrices are well studied in the literature and several constructions and examples are known, see for instance \cite{BKR,K,KS}.

The following is definition of a FSRSG which is introduced in \cite{RAm} and updated in \cite{RA1}.
\noindent  A signed regular graph $\Sigma$ is called \textit{FSRSG}
whenever it satisfies the following conditions:

(i) $|\Sigma|$ is $k$-regular, with $n$ vertices,

(ii) if $\Sigma$ contains at least one positive edge, then there
exist $t\in \mathbb{Z}$ such that $t_{xy}^+-t_{xy}^-=t$, for any edges $xy$,

(iii) there exist $\rho\in \mathbb{Z}$ such that
$\rho_{xy}^+-\rho_{xy}^-=\rho$, for all non-adjacent vertices $x$ and $y$, where $\rho_{xy}^+$, $\rho_{xy}^-$ are the numbers of positive
and negative length two paths joining the vertices. We denote a FSRSG having just mentioned parameters with FSRSG$(n,k,t,\rho)$.

We refer to a signed graph with only two distinct eigenvalues, by \textit{STE}. The following is a general determination of STE's. The result is independently obtained in \cite{RAm} and \cite{ZS}.
\begin{thm}\label{RA1}
  {\rm A signed graph is a STE if and only if it is a FSRSG with $\rho=0$.}
\end{thm}

Recently some problems on the spectrum of signed adjacency matrices have attracted many studies. The authors of \cite{BCKW} has settled few open problems in the spectral theory of signed graphs. The construction problem of STE's is also mentioned there. There  are some results on the problem in literature, see \cite{DM,RA,RAm,ZS}. It is known that the only graphs with two distinct eigenvalues of the ordinary adjacency matrix are the complete graphs, which is extensively far from the known results for the signed graphs. We consider the same problem for signed graphs. We use the well known method, called the \textit{star complement technique}, for construction of STE's. For more details on the notions and definitions see \cite{BR,TAY}. It is well adopted for the family of signed graphs as well as any other symmetric matrices, see \cite{BMS,RAm}. During the preparation of this paper we encounter interesting results from other combinatorial areas such as symmetric weighing matrices, Hadamard matrices, Conference matrices and Ramanujan graphs.
\section{Preliminaries}
The adjacency matrix of a FSRSG satisfies a quadratic equation of the form
\begin{equation}\label{qu1}
A_{\Sigma}^2-tA_{\Sigma}-kI=\rho \overline{A},
\end{equation}
where $\overline{A}$ is the complement of $|\Sigma|$, and $t$, $k$, $\rho$ are as in the definition of FSRSG's. From the equality (\ref{qu1}), for any STE the corresponding $\rho$ parameter leads to be zero, hence for the adjacency matrix of a STE, the following equality must hold.
\begin{equation}\label{qu2}
A_{\Sigma}^2-tA_{\Sigma}-kI=0.
\end{equation}
Therefore the eigenvalues of a STE must be as follows. $$\lambda_1=\frac{t+\sqrt{t^2+4k}}{2},\hspace{0.7 cm}\lambda_2=\frac{t-\sqrt{t^2+4k}}{2}.$$

From now on, we denote $t^2+4k$, by $b$. It is well-known that if a monic integral polynomial admits a non-integral root $\alpha$, then $\alpha$ is irrational and the algebraic conjugate of $\alpha$ is a root of the polynomial with the same multiplicity.  This leads to the following useful lemma.
\begin{lem} \label{eig}
  {\rm Let $\Sigma=(\Gamma,\sigma)$ be a STE, then one of the followings holds.
   \begin{itemize}
     \item $t=0$ and $A_{\Sigma}$ is a weighing matrix of weight $k$, or
     \item $t\neq 0$, then $b$ must be a perfect square.
   \end{itemize}
    }
\end{lem}
\noindent{\bf Proof.} The first part follows by putting $t=0$, in (\ref{qu2}). For the second part as contradiction suppose that $b$ is not a perfect square, so $\frac{-t+\sqrt{b}}{2}$ is irrational. Therefore since the characteristic polynomial of $\Sigma$ is a monic integral polynomial, the algebraic conjugate of $\frac{-t+\sqrt{b}}{2}$, that is, $\frac{t-\sqrt{b}}{2}$ must be an eigenvalue of $\Sigma$ with the same multiplicity, say $m$. By the assumption $\Sigma$ has the spectrum $(\frac{t+\sqrt{b}}{2})^m$, $(\frac{t-\sqrt{b}}{2})^m$. On the other hand trace of $A_{\Sigma}$ is zero, thus we have the following equality,  $$m\frac{t+\sqrt{b}}{2}+m\frac{t-\sqrt{b}}{2}=0,$$ which implies $mt=0$, so $t=0$, a contradiction. Hence $b$ must be a perfect square. $\square$

 \noindent Note that in a signed $k$-regular graphr, $t$ may be an integer in the interval $[-k+1,k-1]$. Regarding the previous statements a triple $(t,\lambda_1,\lambda_2)$ is called \textit{admissible} for a STE, if for $k=-\lambda_1\lambda_2$, either $b=t^2+4k$ is a perfect square or $t$ is zero. With a simple calculation we observe the following.

\begin{thm} \label{para}{\rm The following parameters are all possible admissible $(t,\lambda_1,\lambda_2)$'s of the STE's with valencies $5,6,\ldots,10$.

{\hspace{1.5cm} \begin{tabular}{|c|l | }
   \hline
   $k=5$ & ($4,$ $5,$ $-1$),   ($0,$ $\sqrt{5},$$-\sqrt{5}$),  ($-4,$ $1$,$-5$)  \\ \hline
   $k=6$ & ($5,$  $6$, $-1$),  ($1,$  $3$, $-2$),  ($0,$ $\sqrt{6}$, $-\sqrt{6}$),   ($-1,$   $2$, $-3$),  ($-5,$  $1$, $-6$) \\ \hline
   $k=7$ &  ($6,$ $7,$ $-1$),  ($0,$ $\sqrt{7},$$-\sqrt{7}$),  ($-6,$ $1,$ $-7$)   \\ \hline
   $k=8$ & ($7,$ $8,$ $-1$),  ($2,$ $4,$ $-2$), ($0,$ $\sqrt{8},$$-\sqrt{8}$),  ($-2,$ $-4,$ $2$), ($-7,$ $1,$ $-8$)    \\ \hline
   $k=9$  & ($8,$ $9,$ $-1$),  ($0$, $3$, $-3$),  ($-8,$ $1,$ $-9$)     \\ \hline
   $k=10$ & ($9,$ $10,$ $-1$),  ($3,$ $5,$ $-2$),  ($0,$ $\sqrt{10},$$-\sqrt{10}$),  ($-3,$ $2,$ $-5$),  ($-9,$ $1,$ $-10$)  \\ \hline
  \end{tabular}}
  $$\textbf{Table 1.} \textrm{The admissible parameters for the  STE's with $\{5,6,\ldots,10\}$-regular grounds}$$
}\end{thm}

The notion of line of signed graphs is defined in \cite{ZASee}. The \textit{line graph} of $\Sigma$ is the signed graph $\Lambda(\Sigma)=(\Lambda(G),\sigma_{\Lambda})$, where $\Lambda(G)$ is the ordinary
line graph of the underlying graph and $\sigma_{\Lambda}$ is a sign function such that every cycle in $\Sigma$ becomes a cycle
with the same sign in the line graph, and any three edges incident with a common vertex become a
negative triangle in the line graph. The adjacency matrix is $A_{\Lambda(\Sigma)} = 2I_m-H(\Sigma)^tH(\Sigma)$. Where $H(\Sigma)$ is the \textit{incidence matrix} of a signed graph $\Sigma$ on $n$ vertices with $m$ edges, which is defined in \cite{ZASe}, as follows. It is the $n\times m$ matrix $H(\Sigma)=[\eta_{i,j}]$, where
$$\eta_{i,j}=\left\{
               \begin{array}{ll}
                 0, & \hbox{if $v_i$ is not incident with $e_j$;} \\
                 \pm 1, & \hbox{if $v_i$ is incident with $e_j$.}
               \end{array}
             \right.
$$
such that for an edge $v_iv_j$, the equality $\eta_{ik}\eta_{jk}=\sigma(v_iv_j)$ holds. The incidence matrix
is unique up to negation of columns.
The following theorem from \cite[Theorem 3.10]{GESAZA} will be useful.
\begin{thm} \label{eiglin}
{\rm Assume $\Sigma$ is a signed $k$-regular graph. Let the eigenvalues of $\Sigma$
be $\lambda_1$, $\ldots$, $\lambda_{b(-\Sigma)}=-k$, $-k<\lambda_{b(-\Sigma)+1},\ldots,\lambda_{n-b(\Sigma)}<k$, and $\lambda_{n-b(\Sigma)},\ldots \lambda_n=k$.
The line graph $\Lambda(\Sigma)$ has eigenvalues $\lambda_{b(-\Sigma)+1}-k+2,\ldots, \lambda_{n-b(\Sigma)}-k+2$  and eigenvalue $2$ with multiplicity $m-n+b(\Sigma)$.}
\end{thm}
Note that $b(\Sigma)$ is the number of balanced components of $\Sigma$. 
\subsection{The star complement technique}
In the following we will summarize the method which is extensively surveyed in \cite{BR,TAY}.

Let $A$ be an $n\times n$ matrix, and $\mu$ be an eigenvalue of $A$ having multiplicity $k$. A $k$-subset $X$ of the vertex set is called a \textit{star set} if the sub-matrix of $A$, obtained by removing rows and columns corresponding to $X$ does not have $\mu$ as an eigenvalue. In signed graph context, a star set for an eigenvalue $\mu$ in a signed graph $\Sigma$ is a subset $X$ of vertices such that $\mu$ is not an eigenvalue of $\Sigma\setminus X$. The signed graph $\Sigma\setminus X$ is called a \textit{star complement} for $\mu$ in $\Sigma$. The following theorem, known as the \textit{reconstruction }theorem, will be useful later.
\begin{thm} {\bf \cite{BR}} \label{st} {\rm Let $X$ be a set of $k$ vertices in a signed graph $\Sigma$ and suppose the matrix $A_{\Sigma}$ is of the following form
$$A_{\Sigma}=\left(
  \begin{array}{cc}
    A_X & B \\
    B^t & C \\
  \end{array}
\right),$$ where $A_X$ is the adjacency matrix of the signed graph of $\Sigma$ induced on $X$. Then $X$ is a star set for $\mu$ in $\Sigma$ if and only if $\mu$ is not an eigenvalue of $C$ and $$\mu I-A_X=B(\mu I-C)^{-1}B^t.$$
}
\end{thm}

The star complement method will also be used for the construction of STE's, we state the following immediate consequence of Theorem \ref{st} and the fact that any real symmetric matrix has a \textit{partition to star complements}, \cite{CRS}, with out proof.
\begin{prop} \label{cons}
  {\rm A signed graph $\Sigma$ on $n$ vertices, have spectrum $[\lambda_1^a,\lambda_2^{n-a}]$ if and only if the vertex set of $\Sigma$ has a partition $X,Y$  such that the matrix $A_{\Sigma}$ can be written as the following $$\left(
                                                           \begin{array}{cc}
                                                             A_X & B \\
                                                             B^t & A_Y \\
                                                           \end{array}
                                                         \right),
$$ such that $A_X$ ($A_Y$) doesn't have eigenvalue $\lambda_1$ ($\lambda_2$) and the following equalities hold. $$\lambda_2 I-A_X=B(\lambda_2 I-A_Y)^{-1}B^t, \hspace{1cm}\lambda_1 I-A_Y=B^t(\lambda_1 I-A_X)^{-1}B.$$}
\end{prop}
\section{Constructions}

We have the following consequence of Theorem \ref{eiglin} for the line graph of the all positive complete graph.
\begin{prop} \label{kam}
  {\rm Let $\Sigma$ be the all positive signed graph on the complete graph $K_n$, that is $K_n^+$. Then $\Lambda(K_n^+)$ is a signed graph with the following spectrum. $$[(2-n)^{n-1},2^{{n\choose 2}-n+1}].$$ Moreover the spectrum of the signed graph $-\Lambda(K_n^+)$ is $$[(n-2)^{n-1},(-2)^{{n\choose 2}-n+1}].$$  }
\end{prop}
{\bf Proof.} It follows by Theorem \ref{eiglin} and the fact that $K_n^+$ has spectrum $[(n-1)^1, -1^{n-1}]$. $\square$

Note that the admissible $(t,\lambda_1,\lambda_2)$ parameters of STE's on $k$-regular graphs, with $5\leq k\leq 10$, presented in Table $1$, have the following general types. We will treat each of them seperately.
\begin{itemize}
  \item {\bf Type 1.} $(n-2,n-1,-1)$, or its negation $(-n+2,1,-n+1)$,
  \item {\bf Type 2.} $(\frac{k}{2}-2,\frac{k}{2},-2)$, or its negation $(-\frac{k}{2}+2,2,-\frac{k}{2})$,
  \item {\bf Type 3.} $(0,\sqrt{k},-\sqrt{k})$.
\end{itemize}
Where $n$, $k$ respectively denote the number of vertices, and valency of the signed graph $\Sigma$ which is an STE with the given parameters.

We recall that in a STE's we have $k=-\lambda_1\lambda_2$ and $t=\lambda_1+\lambda_2$. Also from the triples $(t,\lambda_1,\lambda_2)$ and $(-t,-\lambda_1,-\lambda_2)$, we only need to treat one of them, as if $\Sigma$ is a STE with parameters $(t,\lambda_1,\lambda_2)$, then $-\Sigma$ have parameters $(-t,-\lambda_2,-\lambda_1)$. Note also the only signed graphs posing the parameters $(n-2,n-1,-1)$ are the all positive complete graph. Thus we only need to consider the two remained cases.

\subsection{STE's with parameters $\mathbf{(\frac{k}{2}-2,\frac{k}{2},-2)}$}
We will present two family of examples of STE's posing the parameters $(\frac{k}{2}-2,\frac{k}{2},-2)$. First we assert the following relation on the multiplicities of the eigenvalues of STE's having the parameters $(\frac{k}{2}-2,\frac{k}{2},-2)$, which follows by considering the trace of $A_{\Sigma}$. $$m(\frac{k}{2})=\frac{4n}{4+k},\hspace{0.7cm} m(-2)=\frac{nk}{4+k}.$$
\begin{lem}{\rm The followings hold.
\begin{itemize}
  \item The signed graph $\Lambda(K_{5})$ is $6$-regular, having spectrum $[3^4, -2^6]$.
  \item The signed graph $\Lambda(K_{6})$ is $8$-regular, having spectrum $[4^5, -2^{10}]$.
  \item The signed graph $\Lambda(K_{7})$ is $10$-regular, having spectrum $[5^6, -2^{15}]$.
\end{itemize}
  }\end{lem}
{\bf Proof.} The assertion follows by Proposition \ref{kam}. $\square$

The following is a signed graph with spectrum $[3^4, -2^6]$. 

{\hspace{5cm}\begin{tikzpicture}[mystyle/.style={draw,shape=circle}]
\def\ngon{10}
\node[regular polygon,regular polygon sides=\ngon,minimum size=5cm] (p) {};
\foreach\x in {1,...,\ngon}{\node[mystyle] (p\x) at (p.corner \x){};}
\foreach\x in {1,...,\numexpr\ngon-1\relax}
\draw (p1) -- (p6)[red];
\draw (p1) -- (p7)[red];
\draw (p1) -- (p9)[red];
\draw (p1) -- (p5)[blue];
\draw (p1) -- (p8)[blue];
\draw (p1) -- (p10)[blue];
\draw (p2) -- (p4)[blue];
\draw (p2) -- (p7)[blue];
\draw (p2) -- (p8)[blue];
\draw (p2) -- (p3)[red];
\draw (p2) -- (p5)[red];
\draw (p2) -- (p9)[red];
\draw (p3) -- (p4)[blue];
\draw (p3) -- (p5)[blue];
\draw (p3) -- (p9)[blue];
\draw (p3) -- (p6)[red];
\draw (p3) -- (p10)[red];
\draw (p4) -- (p6)[red];
\draw (p4) -- (p10)[red];
\draw (p4) -- (p7)[blue];
\draw (p4) -- (p8)[blue];
\draw (p5) -- (p9)[blue];
\draw (p5) -- (p8)[blue];
\draw (p5) -- (p10)[blue];
\draw (p6) -- (p9)[blue];
\draw (p6) -- (p7)[blue];
\draw (p6) -- (p10)[blue];
\draw (p7) -- (p9)[blue];
\draw (p7) -- (p8)[blue];
\draw (p8) -- (p10)[blue];
\end{tikzpicture}}
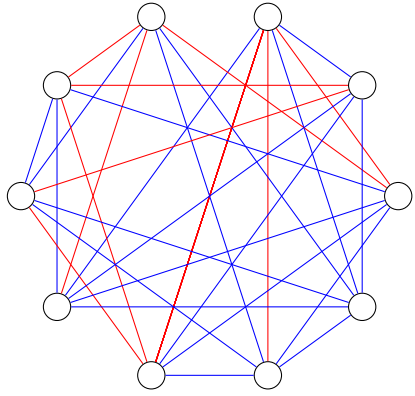
$$\textbf{Figure 1.} \textrm{ The signed graph $-\Lambda(K_5)$}$$

\subsection{$8$-regular signed graphs with spectrum $4^m,-2^{2m}$}

Using Proposition \ref{kam}, we find out the signed graph $\Lambda(K_6)$ has spectrum $4^5,-2^{10}$. At this part using Proposition \ref{cons}, we construct STE's with spectrum $4^m,-2^{2m}$. The following lemma has crucial role in that regard.

\begin{lem} \label{con}
  {\rm Let $W_1$, $W_2$ be two weighing matrices of order $m$ and weight $4$. The following matrix has spectrum $[4^m,-2^{2m}]$.
$$\mathcal{A}(W_1,W_2)=\left(
  \begin{array}{ccc}
    O_m & W_1 & W_2 \\
    W_1^t & O_m & \frac{1}{2}W_1^tW_2 \\
    W_2^t & \frac{1}{2}W_2^tW_1 & O_m \\
  \end{array}
\right).$$
}
\end{lem}
{\bf Proof. }We prove that the subset $\{1,2,\ldots,m\}$ is a star set for $4$, and the subset $\{m+1,m+2,\ldots,3m\}$ is a star set for $-2$. Let $C_1$, $C_2$ be the submatrices corresponding to the mentioned vertex sets respectively. That is, $$C_1=O_m,\textrm{ and }C_2=\left(
                                                                                                               \begin{array}{cc}
                                                                                                                 O_m& \frac{1}{2}W_1^tW_2 \\
                                                                                                                 \frac{1}{2}W_2^tW_1 & O_m\\
                                                                                                               \end{array}
                                                                                                             \right).
$$ We set $B=(W_1|W_2)$, now by Theorem \ref{st}, it suffices to prove the following equalities. $$4I_m=B(4I_{2m}-C_2)^{-1}B^t,$$
$$-2I_{2m}-C_2=B^t(-2I_{m}-C_1)^{-1}B.$$By a simple calculation we obtain the following equality, $$(4I_{2m}-C_2)^{-1}=\left(
                                                                                                                         \begin{array}{cc}
                                                                                                                           \frac{1}{3}I_m & \frac{W_1^tW_2}{24} \\
                                                                                                                           \frac{W_2^tW_1}{24}  & \frac{1}{3}I_m  \\
                                                                                                                         \end{array}
                                                                                                                       \right).
$$
Therefore we obtain the followings by the view of the equalities $W_1W_1^t=W_2W_2^t=4I_m$. $$B(4I_{2m}-C_2)^{-1}B^t=(W_1|W_2)\left(
                                                                                                                         \begin{array}{cc}
                                                                                                                           \frac{1}{3}I_m & \frac{W_1^tW_2}{24} \\
                                                                                                                           \frac{W_2^tW_1}{24}  & \frac{1}{3}I_m  \\
                                                                                                                         \end{array}
                                                                                                                       \right)\left(
                                                                                                                                \begin{array}{c}
                                                                                                                                  W_1^t \\
                                                                                                                                  W_2^t \\
                                                                                                                                \end{array}
                                                                                                                              \right)=(W_1|W_2)\left(\begin{array}{c}
                                                                                                                                  \frac{1}{2}W_1^t \\
                                                                                                                                  \frac{1}{2}W_2^t \\

                                                                                                                                \end{array} \right)=4I_m.$$
Thus the first equality holds, now we prove the second equality. Note that $(-2I_{n}-C_1)^{-1}=\frac{-1}{2}I_n$, therefore we need to prove the following so that the assertion follows. $$B^t(-2I_{m}-C_1)^{-1}B=\frac{-1}{2}\left(
                                                                                                                                \begin{array}{c}
                                                                                                                                  W_1^t \\
                                                                                                                                  W_2^t \\
                                                                                                                                \end{array}
                                                                                                                              \right)(W_1|W_2)=\frac{-1}{2}\left(
                                                                                                                                                 \begin{array}{cc}
                                                                                                                                                   W_1^tW_1 & W_1^tW_2\\
                                                                                                                                                   W_2^tW_1 &  W_2^tW_2  \\
                                                                                                                                                 \end{array}
                                                                                                                                               \right)=-2I_{2m}-C_2.
$$
Now the assertion follows by Theorem \ref{st}. $\square$

\begin{cor}
{\rm If the weighing matrices $W_1$ and $W_2$ are chosen so that the matrix $\frac{1}{2}W_1^tW_2$ has $0,\pm 1$ entries, or equivalently the matrix $W_1^tW_2$ has $0,\pm 2$ entries, then the matrix $\mathcal{A}(W_1,W_2)$ is corresponding to a $8$-regular graph on $3m$ vertices, having spectrum $[4^m,-2^{2m}]$.}
\end{cor}

If two weighing matrices $W_1$ and $W_2$ of weight $4$ are so that the matrix $W_1^tW_2$ has $0,\pm 2$ entries, then we call them \textit{semi-orthogonal} weighing matrices. In the following we provide examples of semi-orthogonal weighing matrices of any even order $m\geq 6$. We should mention that a family of weighing matrices of even order is presented in \cite{BKR}, but our method provides a different matrix. We make use of the following two pattern matrices $X$ and $Y$ of order $\frac{m}{2}$. For $i,j=1,2,\ldots,\frac{m}{2},$ the $(i,j)$'th entries of $X$ and $Y$ follows.
$$
    \begin{array}{cc}
      X(i,j)=\left\{
           \begin{array}{ll}
             1, &j-i\equiv 0,1 \hbox{ mod }\frac{m}{2};  \\
             0, & \hbox{other wise.}
           \end{array}
         \right. & \hspace{0.5cm}Y(i,j)=\left\{
           \begin{array}{ll}
             1, &j-i\equiv 1, \frac{m}{2}-2 \hbox{ mod }\frac{m}{2};  \\
             0, & \hbox{other wise.}
           \end{array}
         \right. \\
    \end{array}
$$
Note that the matrices $X,Y$ have two entry $1$ at each row and each column. Moreover $X,Y$ do not share a common entry $1$. We define the matrix $W_1$, similarly $W_2$ as the following.
\begin{itemize}
 \item The first $\frac{m}{2}\times m$ block of $W_1$, say $F=F_X$, is filled as follows. In fact the matrix $X$ is expanded to the matrix $F$ of size $\frac{m}{2}\times m$. For $i,j=1,2,\ldots,\frac{m}{2},$ if $X(i,j)=0$, then replace the $(i,j)$'th entry of $X$ with $[0,0]$, for $j=1,2,\ldots,m/2$ if in the $j$'th column the one entries occur in the coordinates $(i_0,j)$, and $(i_1,j)$, where $i_0<i_1$, then replace the $(i_0,j)$'th entry with $[1,1]$, and the $(i_1,j)$'th entries with $[1,-1]$.
 \item The next $\frac{m}{2}\times m$ block of $W_1$, say $R=R_X$, is filled as follows. The entries of $R_X$ in the corresponding zero coordinates of $F_X$ are zero. For $i=1,2,\ldots,\frac{m}{2},$ if non-zero entries of $F_X$ occur in the coordinates $(i,j_0)$, $(i,j_1)$, $(i,j_2)$, $(i,j_3)$, where $j_0<j_1<j_2<j_3$, then we set $$R_X(i,j_0)=F_X(i,j_0), R_X(i,j_1)=F_X(i,j_1),$$$$ R_X(i,j_2)=-F_X(i,j_2), R_X(i,j_3)=-F_X(i,j_3).$$
 
 \end{itemize}
Now set $W_1$ to be the $m\times m$ matrix, $\left(
                       \begin{array}{c}
                         F \\
                         R \\
                       \end{array}
                     \right)
$. The matrix $W_2$ is constructed similarly by considering $Y$ instead of $X$. See the following matrices for illustration.
$${\scriptsize X=\left(
      \begin{array}{ccccccc}
        1 &\hspace{-0.2cm} 0 &\hspace{-0.2cm} 0 &\hspace{-0.2cm} 0 &\hspace{-0.2cm} 0 &\hspace{-0.2cm} 0 &\hspace{-0.2cm} 1 \\
        1 &\hspace{-0.2cm} 1 &\hspace{-0.2cm} 0 &\hspace{-0.2cm} 0 &\hspace{-0.2cm} 0 &\hspace{-0.2cm} 0 &\hspace{-0.2cm} 0 \\
        0 &\hspace{-0.2cm} 1 &\hspace{-0.2cm} 1 &\hspace{-0.2cm} 0 &\hspace{-0.2cm} 0 &\hspace{-0.2cm} 0 &\hspace{-0.2cm} 0 \\
        0 &\hspace{-0.2cm} 0 &\hspace{-0.2cm} 1 &\hspace{-0.2cm} 1 &\hspace{-0.2cm} 0 &\hspace{-0.2cm} 0 &\hspace{-0.2cm} 0 \\
        0 &\hspace{-0.2cm} 0 &\hspace{-0.2cm} 0 &\hspace{-0.2cm} 1 &\hspace{-0.2cm} 1 &\hspace{-0.2cm} 0 &\hspace{-0.2cm} 0 \\
        0 &\hspace{-0.2cm} 0 &\hspace{-0.2cm} 0 &\hspace{-0.2cm} 0 &\hspace{-0.2cm} 1 &\hspace{-0.2cm} 1 &\hspace{-0.2cm} 0 \\
        0 &\hspace{-0.2cm} 0 &\hspace{-0.2cm} 0 &\hspace{-0.2cm} 0 &\hspace{-0.2cm} 0 &\hspace{-0.2cm} 1 &\hspace{-0.2cm} 1 \\
      \end{array}
    \right),\hspace{1cm}Y=\left(
      \begin{array}{ccccccc}
        0 &\hspace{-0.2cm} 1 &\hspace{-0.2cm} 0 &\hspace{-0.2cm} 0 &\hspace{-0.2cm} 0 &\hspace{-0.2cm} 1 &\hspace{-0.2cm} 0 \\
        0 &\hspace{-0.2cm} 0 &\hspace{-0.2cm} 1 &\hspace{-0.2cm} 0 &\hspace{-0.2cm} 0 &\hspace{-0.2cm} 0 &\hspace{-0.2cm} 1 \\
        1 &\hspace{-0.2cm} 0 &\hspace{-0.2cm} 0 &\hspace{-0.2cm} 1 &\hspace{-0.2cm} 0 &\hspace{-0.2cm} 0 &\hspace{-0.2cm} 0 \\
        0 &\hspace{-0.2cm} 1 &\hspace{-0.2cm} 0 &\hspace{-0.2cm} 0 &\hspace{-0.2cm} 1 &\hspace{-0.2cm} 0 &\hspace{-0.2cm} 0 \\
        0 &\hspace{-0.2cm} 0 &\hspace{-0.2cm} 1 &\hspace{-0.2cm} 0 &\hspace{-0.2cm} 0 &\hspace{-0.2cm} 1 &\hspace{-0.2cm} 0 \\
        0 &\hspace{-0.2cm} 0 &\hspace{-0.2cm} 0 &\hspace{-0.2cm} 1 &\hspace{-0.2cm} 0 &\hspace{-0.2cm} 0 &\hspace{-0.2cm} 1 \\
        1 &\hspace{-0.2cm} 0 &\hspace{-0.2cm} 0 &\hspace{-0.2cm} 0 &\hspace{-0.2cm} 1 &\hspace{-0.2cm} 0 &\hspace{-0.2cm} 0 \\
      \end{array}
    \right).}
$$
$${\scriptsize F_X=\left(
        \begin{array}{cccccccccccccc}
          1 &\hspace{-0.2cm} 1 &\hspace{-0.2cm} 0 &\hspace{-0.2cm} 0 &\hspace{-0.2cm} 0 &\hspace{-0.2cm} 0 &\hspace{-0.2cm} 0 &\hspace{-0.2cm} 0 &\hspace{-0.2cm} 0 &\hspace{-0.2cm} 0 &\hspace{-0.2cm} 0 &\hspace{-0.2cm} 0 &\hspace{-0.2cm} 1 &\hspace{-0.2cm} 1 \\
          1 &\hspace{-0.2cm} - &\hspace{-0.2cm} 1 &\hspace{-0.2cm} 1 &\hspace{-0.2cm} 0 &\hspace{-0.2cm} 0 &\hspace{-0.2cm} 0 &\hspace{-0.2cm} 0 &\hspace{-0.2cm} 0 &\hspace{-0.2cm} 0 &\hspace{-0.2cm} 0 &\hspace{-0.2cm} 0 &\hspace{-0.2cm} 0 &\hspace{-0.2cm} 0 \\
          0 &\hspace{-0.2cm} 0 &\hspace{-0.2cm} 1 &\hspace{-0.2cm} - &\hspace{-0.2cm} 1 &\hspace{-0.2cm} 1 &\hspace{-0.2cm} 0 &\hspace{-0.2cm} 0 &\hspace{-0.2cm} 0 &\hspace{-0.2cm} 0 &\hspace{-0.2cm} 0 &\hspace{-0.2cm} 0 &\hspace{-0.2cm} 0 &\hspace{-0.2cm} 0 \\
          0 &\hspace{-0.2cm} 0 &\hspace{-0.2cm} 0 &\hspace{-0.2cm} 0 &\hspace{-0.2cm} 1 &\hspace{-0.2cm} - &\hspace{-0.2cm} 1 &\hspace{-0.2cm} 1 &\hspace{-0.2cm} 0 &\hspace{-0.2cm} 0 &\hspace{-0.2cm} 0 &\hspace{-0.2cm} 0 &\hspace{-0.2cm} 0 &\hspace{-0.2cm} 0 \\
          0 &\hspace{-0.2cm} 0 &\hspace{-0.2cm} 0 &\hspace{-0.2cm} 0 &\hspace{-0.2cm} 0 &\hspace{-0.2cm} 0 &\hspace{-0.2cm} 1 &\hspace{-0.2cm} - &\hspace{-0.2cm} 1 &\hspace{-0.2cm} 1 &\hspace{-0.2cm} 0 &\hspace{-0.2cm} 0 &\hspace{-0.2cm} 0 &\hspace{-0.2cm} 0 \\
          0 &\hspace{-0.2cm} 0 &\hspace{-0.2cm} 0 &\hspace{-0.2cm} 0 &\hspace{-0.2cm} 0 &\hspace{-0.2cm} 0 &\hspace{-0.2cm} 0 &\hspace{-0.2cm} 0 &\hspace{-0.2cm} 1 &\hspace{-0.2cm} - &\hspace{-0.2cm} 1 &\hspace{-0.2cm} 1 &\hspace{-0.2cm} 0 &\hspace{-0.2cm} 0 \\
          0 &\hspace{-0.2cm} 0 &\hspace{-0.2cm} 0 &\hspace{-0.2cm} 0 &\hspace{-0.2cm} 0 &\hspace{-0.2cm} 0 &\hspace{-0.2cm} 0 &\hspace{-0.2cm} 0 &\hspace{-0.2cm} 0 &\hspace{-0.2cm} 0 &\hspace{-0.2cm} 1 &\hspace{-0.2cm} - &\hspace{-0.2cm} 1 &\hspace{-0.2cm} - \\
        \end{array}
      \right)\hspace{0.5cm}R_X=\left(
        \begin{array}{cccccccccccccc}
          1 &\hspace{-0.2cm} 1 &\hspace{-0.2cm} 0 &\hspace{-0.2cm} 0 &\hspace{-0.2cm} 0 &\hspace{-0.2cm} 0 &\hspace{-0.2cm} 0 &\hspace{-0.2cm} 0 &\hspace{-0.2cm} 0 &\hspace{-0.2cm} 0 &\hspace{-0.2cm} 0 &\hspace{-0.2cm} 0 &\hspace{-0.2cm} - &\hspace{-0.2cm} - \\
          1 &\hspace{-0.2cm} - &\hspace{-0.2cm} - &\hspace{-0.2cm} - &\hspace{-0.2cm} 0 &\hspace{-0.2cm} 0 &\hspace{-0.2cm} 0 &\hspace{-0.2cm} 0 &\hspace{-0.2cm} 0 &\hspace{-0.2cm} 0 &\hspace{-0.2cm} 0 &\hspace{-0.2cm} 0 &\hspace{-0.2cm} 0 &\hspace{-0.2cm} 0 \\
          0 &\hspace{-0.2cm} 0 &\hspace{-0.2cm} 1 &\hspace{-0.2cm} - &\hspace{-0.2cm} - &\hspace{-0.2cm} - &\hspace{-0.2cm} 0 &\hspace{-0.2cm} 0 &\hspace{-0.2cm} 0 &\hspace{-0.2cm} 0 &\hspace{-0.2cm} 0 &\hspace{-0.2cm} 0 &\hspace{-0.2cm} 0 &\hspace{-0.2cm} 0 \\
          0 &\hspace{-0.2cm} 0 &\hspace{-0.2cm} 0 &\hspace{-0.2cm} 0 &\hspace{-0.2cm} 1 &\hspace{-0.2cm} - &\hspace{-0.2cm} - &\hspace{-0.2cm} - &\hspace{-0.2cm} 0 &\hspace{-0.2cm} 0 &\hspace{-0.2cm} 0 &\hspace{-0.2cm} 0 &\hspace{-0.2cm} 0 &\hspace{-0.2cm} 0 \\
          0 &\hspace{-0.2cm} 0 &\hspace{-0.2cm} 0 &\hspace{-0.2cm} 0 &\hspace{-0.2cm} 0 &\hspace{-0.2cm} 0 &\hspace{-0.2cm} 1 &\hspace{-0.2cm} - &\hspace{-0.2cm} - &\hspace{-0.2cm} - &\hspace{-0.2cm} 0 &\hspace{-0.2cm} 0 &\hspace{-0.2cm} 0 &\hspace{-0.2cm} 0 \\
          0 &\hspace{-0.2cm} 0 &\hspace{-0.2cm} 0 &\hspace{-0.2cm} 0 &\hspace{-0.2cm} 0 &\hspace{-0.2cm} 0 &\hspace{-0.2cm} 0 &\hspace{-0.2cm} 0 &\hspace{-0.2cm} 1 &\hspace{-0.2cm} - &\hspace{-0.2cm} - &\hspace{-0.2cm} - &\hspace{-0.2cm} 0 &\hspace{-0.2cm} 0 \\
          0 &\hspace{-0.2cm} 0 &\hspace{-0.2cm} 0 &\hspace{-0.2cm} 0 &\hspace{-0.2cm} 0 &\hspace{-0.2cm} 0 &\hspace{-0.2cm} 0 &\hspace{-0.2cm} 0 &\hspace{-0.2cm} 0 &\hspace{-0.2cm} 0 &\hspace{-0.2cm} 1 &\hspace{-0.2cm} - &\hspace{-0.2cm} - &\hspace{-0.2cm} 1 \\
        \end{array}
      \right)}
$$

$${\scriptsize F_Y=\left(
      \begin{array}{cccccccccccccc}
        0&\hspace{-0.2cm}0 &\hspace{-0.2cm} 1&\hspace{-0.2cm}1 &\hspace{-0.2cm} 0 &\hspace{-0.2cm} 0 &\hspace{-0.2cm} 0&\hspace{-0.2cm} 0 &\hspace{-0.2cm} 0 &\hspace{-0.2cm} 0 &\hspace{-0.2cm} 1&\hspace{-0.2cm}1 &\hspace{-0.2cm} 0&\hspace{-0.2cm}0 \\
        0 &\hspace{-0.2cm} 0 &\hspace{-0.2cm}0 &\hspace{-0.2cm} 0 &\hspace{-0.2cm} 1&\hspace{-0.2cm}1 &\hspace{-0.2cm} 0 &\hspace{-0.2cm} 0 &\hspace{-0.2cm} 0&\hspace{-0.2cm} 0 &\hspace{-0.2cm} 0 &\hspace{-0.2cm} 0 &\hspace{-0.2cm} 1&\hspace{-0.2cm} 1 \\
       1 &\hspace{-0.2cm} 1 &\hspace{-0.2cm} 0 &\hspace{-0.2cm} 0&\hspace{-0.2cm} 0 &\hspace{-0.2cm} 0 &\hspace{-0.2cm} 1&\hspace{-0.2cm} 1 &\hspace{-0.2cm} 0 &\hspace{-0.2cm} 0 &\hspace{-0.2cm} 0&\hspace{-0.2cm} 0 &\hspace{-0.2cm} 0 &\hspace{-0.2cm} 0 \\
        0 &\hspace{-0.2cm}0 &\hspace{-0.2cm} 1&\hspace{-0.2cm}- &\hspace{-0.2cm} 0 &\hspace{-0.2cm} 0&\hspace{-0.2cm} 0 &\hspace{-0.2cm} 0 &\hspace{-0.2cm} 1&\hspace{-0.2cm} 1 &\hspace{-0.2cm} 0 &\hspace{-0.2cm} 0&\hspace{-0.2cm} 0 &\hspace{-0.2cm} 0 \\
        0 &\hspace{-0.2cm} 0 &\hspace{-0.2cm}0 &\hspace{-0.2cm} 0 &\hspace{-0.2cm} 1&\hspace{-0.2cm}- &\hspace{-0.2cm} 0 &\hspace{-0.2cm} 0&\hspace{-0.2cm} 0 &\hspace{-0.2cm} 0 &\hspace{-0.2cm} 1&\hspace{-0.2cm}- &\hspace{-0.2cm} 0&\hspace{-0.2cm} 0 \\
        0 &\hspace{-0.2cm} 0 &\hspace{-0.2cm} 0 &\hspace{-0.2cm}0 &\hspace{-0.2cm} 0 &\hspace{-0.2cm} 0 &\hspace{-0.2cm} 1&\hspace{-0.2cm}- &\hspace{-0.2cm} 0 &\hspace{-0.2cm} 0&\hspace{-0.2cm} 0 &\hspace{-0.2cm} 0 &\hspace{-0.2cm} 1&\hspace{-0.2cm}- \\
        1&\hspace{-0.2cm}- &\hspace{-0.2cm} 0 &\hspace{-0.2cm} 0 &\hspace{-0.2cm} 0 &\hspace{-0.2cm} 0 &\hspace{-0.2cm} 0 &\hspace{-0.2cm} 0 &\hspace{-0.2cm} 1&\hspace{-0.2cm}- &\hspace{-0.2cm} 0 &\hspace{-0.2cm} 0&\hspace{-0.2cm} 0 &\hspace{-0.2cm} 0 \\
      \end{array}
    \right), \hspace{0.5cm}R_Y=\left(
      \begin{array}{cccccccccccccc}
        0&\hspace{-0.2cm}0 &\hspace{-0.2cm} 1&\hspace{-0.2cm}1 &\hspace{-0.2cm} 0 &\hspace{-0.2cm} 0 &\hspace{-0.2cm} 0&\hspace{-0.2cm} 0 &\hspace{-0.2cm} 0 &\hspace{-0.2cm} 0 &\hspace{-0.2cm} -&\hspace{-0.2cm}- &\hspace{-0.2cm} 0&\hspace{-0.2cm}0 \\
        0 &\hspace{-0.2cm} 0 &\hspace{-0.2cm}0 &\hspace{-0.2cm} 0 &\hspace{-0.2cm} 1&\hspace{-0.2cm}1 &\hspace{-0.2cm} 0 &\hspace{-0.2cm} 0 &\hspace{-0.2cm} 0&\hspace{-0.2cm} 0 &\hspace{-0.2cm} 0 &\hspace{-0.2cm} 0 &\hspace{-0.2cm} -&\hspace{-0.2cm} - \\
       1 &\hspace{-0.2cm} 1 &\hspace{-0.2cm} 0 &\hspace{-0.2cm} 0&\hspace{-0.2cm} 0 &\hspace{-0.2cm} 0 &\hspace{-0.2cm} -&\hspace{-0.2cm} - &\hspace{-0.2cm} 0 &\hspace{-0.2cm} 0 &\hspace{-0.2cm} 0&\hspace{-0.2cm} 0 &\hspace{-0.2cm} 0 &\hspace{-0.2cm} 0 \\
        0 &\hspace{-0.2cm}0 &\hspace{-0.2cm} 1&\hspace{-0.2cm}- &\hspace{-0.2cm} 0 &\hspace{-0.2cm} 0&\hspace{-0.2cm} 0 &\hspace{-0.2cm} 0 &\hspace{-0.2cm} -&\hspace{-0.2cm} - &\hspace{-0.2cm} 0 &\hspace{-0.2cm} 0&\hspace{-0.2cm} 0 &\hspace{-0.2cm} 0 \\
        0 &\hspace{-0.2cm} 0 &\hspace{-0.2cm}0 &\hspace{-0.2cm} 0 &\hspace{-0.2cm} 1&\hspace{-0.2cm}- &\hspace{-0.2cm} 0 &\hspace{-0.2cm} 0&\hspace{-0.2cm} 0 &\hspace{-0.2cm} 0 &\hspace{-0.2cm} -&\hspace{-0.2cm}1 &\hspace{-0.2cm} 0&\hspace{-0.2cm} 0 \\
        0 &\hspace{-0.2cm} 0 &\hspace{-0.2cm} 0 &\hspace{-0.2cm}0 &\hspace{-0.2cm} 0 &\hspace{-0.2cm} 0 &\hspace{-0.2cm} 1&\hspace{-0.2cm}- &\hspace{-0.2cm} 0 &\hspace{-0.2cm} 0&\hspace{-0.2cm} 0 &\hspace{-0.2cm} 0 &\hspace{-0.2cm} -&\hspace{-0.2cm}1 \\
        1&\hspace{-0.2cm}- &\hspace{-0.2cm} 0 &\hspace{-0.2cm} 0 &\hspace{-0.2cm} 0 &\hspace{-0.2cm} 0 &\hspace{-0.2cm} 0 &\hspace{-0.2cm} 0 &\hspace{-0.2cm} -&\hspace{-0.2cm}1 &\hspace{-0.2cm} 0 &\hspace{-0.2cm} 0&\hspace{-0.2cm} 0 &\hspace{-0.2cm} 0 \\
      \end{array}
    \right).}$$
$${\scriptsize W_1=\left(
        \begin{array}{cccccccccccccc}
          1 &\hspace{-0.2cm} 1 &\hspace{-0.2cm} 0 &\hspace{-0.2cm} 0 &\hspace{-0.2cm} 0 &\hspace{-0.2cm} 0 &\hspace{-0.2cm} 0 &\hspace{-0.2cm} 0 &\hspace{-0.2cm} 0 &\hspace{-0.2cm} 0 &\hspace{-0.2cm} 0 &\hspace{-0.2cm} 0 &\hspace{-0.2cm} 1 &\hspace{-0.2cm} 1 \\
          1 &\hspace{-0.2cm} - &\hspace{-0.2cm} 1 &\hspace{-0.2cm} 1 &\hspace{-0.2cm} 0 &\hspace{-0.2cm} 0 &\hspace{-0.2cm} 0 &\hspace{-0.2cm} 0 &\hspace{-0.2cm} 0 &\hspace{-0.2cm} 0 &\hspace{-0.2cm} 0 &\hspace{-0.2cm} 0 &\hspace{-0.2cm} 0 &\hspace{-0.2cm} 0 \\
          0 &\hspace{-0.2cm} 0 &\hspace{-0.2cm} 1 &\hspace{-0.2cm} - &\hspace{-0.2cm} 1 &\hspace{-0.2cm} 1 &\hspace{-0.2cm} 0 &\hspace{-0.2cm} 0 &\hspace{-0.2cm} 0 &\hspace{-0.2cm} 0 &\hspace{-0.2cm} 0 &\hspace{-0.2cm} 0 &\hspace{-0.2cm} 0 &\hspace{-0.2cm} 0 \\
          0 &\hspace{-0.2cm} 0 &\hspace{-0.2cm} 0 &\hspace{-0.2cm} 0 &\hspace{-0.2cm} 1 &\hspace{-0.2cm} - &\hspace{-0.2cm} 1 &\hspace{-0.2cm} 1 &\hspace{-0.2cm} 0 &\hspace{-0.2cm} 0 &\hspace{-0.2cm} 0 &\hspace{-0.2cm} 0 &\hspace{-0.2cm} 0 &\hspace{-0.2cm} 0 \\
          0 &\hspace{-0.2cm} 0 &\hspace{-0.2cm} 0 &\hspace{-0.2cm} 0 &\hspace{-0.2cm} 0 &\hspace{-0.2cm} 0 &\hspace{-0.2cm} 1 &\hspace{-0.2cm} - &\hspace{-0.2cm} 1 &\hspace{-0.2cm} 1 &\hspace{-0.2cm} 0 &\hspace{-0.2cm} 0 &\hspace{-0.2cm} 0 &\hspace{-0.2cm} 0 \\
          0 &\hspace{-0.2cm} 0 &\hspace{-0.2cm} 0 &\hspace{-0.2cm} 0 &\hspace{-0.2cm} 0 &\hspace{-0.2cm} 0 &\hspace{-0.2cm} 0 &\hspace{-0.2cm} 0 &\hspace{-0.2cm} 1 &\hspace{-0.2cm} - &\hspace{-0.2cm} 1 &\hspace{-0.2cm} 1 &\hspace{-0.2cm} 0 &\hspace{-0.2cm} 0 \\
          0 &\hspace{-0.2cm} 0 &\hspace{-0.2cm} 0 &\hspace{-0.2cm} 0 &\hspace{-0.2cm} 0 &\hspace{-0.2cm} 0 &\hspace{-0.2cm} 0 &\hspace{-0.2cm} 0 &\hspace{-0.2cm} 0 &\hspace{-0.2cm} 0 &\hspace{-0.2cm} 1 &\hspace{-0.2cm} - &\hspace{-0.2cm} 1 &\hspace{-0.2cm} - \\
          1 &\hspace{-0.2cm} 1 &\hspace{-0.2cm} 0 &\hspace{-0.2cm} 0 &\hspace{-0.2cm} 0 &\hspace{-0.2cm} 0 &\hspace{-0.2cm} 0 &\hspace{-0.2cm} 0 &\hspace{-0.2cm} 0 &\hspace{-0.2cm} 0 &\hspace{-0.2cm} 0 &\hspace{-0.2cm} 0 &\hspace{-0.2cm} - &\hspace{-0.2cm} - \\
          1 &\hspace{-0.2cm} - &\hspace{-0.2cm} - &\hspace{-0.2cm} - &\hspace{-0.2cm} 0 &\hspace{-0.2cm} 0 &\hspace{-0.2cm} 0 &\hspace{-0.2cm} 0 &\hspace{-0.2cm} 0 &\hspace{-0.2cm} 0 &\hspace{-0.2cm} 0 &\hspace{-0.2cm} 0 &\hspace{-0.2cm} 0 &\hspace{-0.2cm} 0 \\
          0 &\hspace{-0.2cm} 0 &\hspace{-0.2cm} 1 &\hspace{-0.2cm} - &\hspace{-0.2cm} - &\hspace{-0.2cm} - &\hspace{-0.2cm} 0 &\hspace{-0.2cm} 0 &\hspace{-0.2cm} 0 &\hspace{-0.2cm} 0 &\hspace{-0.2cm} 0 &\hspace{-0.2cm} 0 &\hspace{-0.2cm} 0 &\hspace{-0.2cm} 0 \\
          0 &\hspace{-0.2cm} 0 &\hspace{-0.2cm} 0 &\hspace{-0.2cm} 0 &\hspace{-0.2cm} 1 &\hspace{-0.2cm} - &\hspace{-0.2cm} - &\hspace{-0.2cm} - &\hspace{-0.2cm} 0 &\hspace{-0.2cm} 0 &\hspace{-0.2cm} 0 &\hspace{-0.2cm} 0 &\hspace{-0.2cm} 0 &\hspace{-0.2cm} 0 \\
          0 &\hspace{-0.2cm} 0 &\hspace{-0.2cm} 0 &\hspace{-0.2cm} 0 &\hspace{-0.2cm} 0 &\hspace{-0.2cm} 0 &\hspace{-0.2cm} 1 &\hspace{-0.2cm} - &\hspace{-0.2cm} - &\hspace{-0.2cm} - &\hspace{-0.2cm} 0 &\hspace{-0.2cm} 0 &\hspace{-0.2cm} 0 &\hspace{-0.2cm} 0 \\
          0 &\hspace{-0.2cm} 0 &\hspace{-0.2cm} 0 &\hspace{-0.2cm} 0 &\hspace{-0.2cm} 0 &\hspace{-0.2cm} 0 &\hspace{-0.2cm} 0 &\hspace{-0.2cm} 0 &\hspace{-0.2cm} 1 &\hspace{-0.2cm} - &\hspace{-0.2cm} - &\hspace{-0.2cm} - &\hspace{-0.2cm} 0 &\hspace{-0.2cm} 0 \\
          0 &\hspace{-0.2cm} 0 &\hspace{-0.2cm} 0 &\hspace{-0.2cm} 0 &\hspace{-0.2cm} 0 &\hspace{-0.2cm} 0 &\hspace{-0.2cm} 0 &\hspace{-0.2cm} 0 &\hspace{-0.2cm} 0 &\hspace{-0.2cm} 0 &\hspace{-0.2cm} 1 &\hspace{-0.2cm} - &\hspace{-0.2cm} - &\hspace{-0.2cm} 1 \\
        \end{array}
      \right),\hspace{0.5cm}W_2=\left(
      \begin{array}{cccccccccccccc}
        0&\hspace{-0.2cm}0 &\hspace{-0.2cm} 1&\hspace{-0.2cm}1 &\hspace{-0.2cm} 0 &\hspace{-0.2cm} 0 &\hspace{-0.2cm} 0&\hspace{-0.2cm} 0 &\hspace{-0.2cm} 0 &\hspace{-0.2cm} 0 &\hspace{-0.2cm} 1&\hspace{-0.2cm}1 &\hspace{-0.2cm} 0&\hspace{-0.2cm}0 \\
        0 &\hspace{-0.2cm} 0 &\hspace{-0.2cm}0 &\hspace{-0.2cm} 0 &\hspace{-0.2cm} 1&\hspace{-0.2cm}1 &\hspace{-0.2cm} 0 &\hspace{-0.2cm} 0 &\hspace{-0.2cm} 0&\hspace{-0.2cm} 0 &\hspace{-0.2cm} 0 &\hspace{-0.2cm} 0 &\hspace{-0.2cm} 1&\hspace{-0.2cm} 1 \\
       1 &\hspace{-0.2cm} 1 &\hspace{-0.2cm} 0 &\hspace{-0.2cm} 0&\hspace{-0.2cm} 0 &\hspace{-0.2cm} 0 &\hspace{-0.2cm} 1&\hspace{-0.2cm} 1 &\hspace{-0.2cm} 0 &\hspace{-0.2cm} 0 &\hspace{-0.2cm} 0&\hspace{-0.2cm} 0 &\hspace{-0.2cm} 0 &\hspace{-0.2cm} 0 \\
        0 &\hspace{-0.2cm}0 &\hspace{-0.2cm} 1&\hspace{-0.2cm}- &\hspace{-0.2cm} 0 &\hspace{-0.2cm} 0&\hspace{-0.2cm} 0 &\hspace{-0.2cm} 0 &\hspace{-0.2cm} 1&\hspace{-0.2cm} 1 &\hspace{-0.2cm} 0 &\hspace{-0.2cm} 0&\hspace{-0.2cm} 0 &\hspace{-0.2cm} 0 \\
        0 &\hspace{-0.2cm} 0 &\hspace{-0.2cm}0 &\hspace{-0.2cm} 0 &\hspace{-0.2cm} 1&\hspace{-0.2cm}- &\hspace{-0.2cm} 0 &\hspace{-0.2cm} 0&\hspace{-0.2cm} 0 &\hspace{-0.2cm} 0 &\hspace{-0.2cm} 1&\hspace{-0.2cm}- &\hspace{-0.2cm} 0&\hspace{-0.2cm} 0 \\
        0 &\hspace{-0.2cm} 0 &\hspace{-0.2cm} 0 &\hspace{-0.2cm}0 &\hspace{-0.2cm} 0 &\hspace{-0.2cm} 0 &\hspace{-0.2cm} 1&\hspace{-0.2cm}- &\hspace{-0.2cm} 0 &\hspace{-0.2cm} 0&\hspace{-0.2cm} 0 &\hspace{-0.2cm} 0 &\hspace{-0.2cm} 1&\hspace{-0.2cm}- \\
        1&\hspace{-0.2cm}- &\hspace{-0.2cm} 0 &\hspace{-0.2cm} 0 &\hspace{-0.2cm} 0 &\hspace{-0.2cm} 0 &\hspace{-0.2cm} 0 &\hspace{-0.2cm} 0 &\hspace{-0.2cm} 1&\hspace{-0.2cm}- &\hspace{-0.2cm} 0 &\hspace{-0.2cm} 0&\hspace{-0.2cm} 0 &\hspace{-0.2cm} 0 \\
        0&\hspace{-0.2cm}0 &\hspace{-0.2cm} 1&\hspace{-0.2cm}1 &\hspace{-0.2cm} 0 &\hspace{-0.2cm} 0 &\hspace{-0.2cm} 0&\hspace{-0.2cm} 0 &\hspace{-0.2cm} 0 &\hspace{-0.2cm} 0 &\hspace{-0.2cm} -&\hspace{-0.2cm}- &\hspace{-0.2cm} 0&\hspace{-0.2cm}0 \\
        0 &\hspace{-0.2cm} 0 &\hspace{-0.2cm}0 &\hspace{-0.2cm} 0 &\hspace{-0.2cm} 1&\hspace{-0.2cm}1 &\hspace{-0.2cm} 0 &\hspace{-0.2cm} 0 &\hspace{-0.2cm} 0&\hspace{-0.2cm} 0 &\hspace{-0.2cm} 0 &\hspace{-0.2cm} 0 &\hspace{-0.2cm} -&\hspace{-0.2cm} - \\
       1 &\hspace{-0.2cm} 1 &\hspace{-0.2cm} 0 &\hspace{-0.2cm} 0&\hspace{-0.2cm} 0 &\hspace{-0.2cm} 0 &\hspace{-0.2cm} -&\hspace{-0.2cm} - &\hspace{-0.2cm} 0 &\hspace{-0.2cm} 0 &\hspace{-0.2cm} 0&\hspace{-0.2cm} 0 &\hspace{-0.2cm} 0 &\hspace{-0.2cm} 0 \\
        0 &\hspace{-0.2cm}0 &\hspace{-0.2cm} 1&\hspace{-0.2cm}- &\hspace{-0.2cm} 0 &\hspace{-0.2cm} 0&\hspace{-0.2cm} 0 &\hspace{-0.2cm} 0 &\hspace{-0.2cm} -&\hspace{-0.2cm} - &\hspace{-0.2cm} 0 &\hspace{-0.2cm} 0&\hspace{-0.2cm} 0 &\hspace{-0.2cm} 0 \\
        0 &\hspace{-0.2cm} 0 &\hspace{-0.2cm}0 &\hspace{-0.2cm} 0 &\hspace{-0.2cm} 1&\hspace{-0.2cm}- &\hspace{-0.2cm} 0 &\hspace{-0.2cm} 0&\hspace{-0.2cm} 0 &\hspace{-0.2cm} 0 &\hspace{-0.2cm} -&\hspace{-0.2cm}1 &\hspace{-0.2cm} 0&\hspace{-0.2cm} 0 \\
        0 &\hspace{-0.2cm} 0 &\hspace{-0.2cm} 0 &\hspace{-0.2cm}0 &\hspace{-0.2cm} 0 &\hspace{-0.2cm} 0 &\hspace{-0.2cm} 1&\hspace{-0.2cm}- &\hspace{-0.2cm} 0 &\hspace{-0.2cm} 0&\hspace{-0.2cm} 0 &\hspace{-0.2cm} 0 &\hspace{-0.2cm} -&\hspace{-0.2cm}1 \\
        1&\hspace{-0.2cm}- &\hspace{-0.2cm} 0 &\hspace{-0.2cm} 0 &\hspace{-0.2cm} 0 &\hspace{-0.2cm} 0 &\hspace{-0.2cm} 0 &\hspace{-0.2cm} 0 &\hspace{-0.2cm} -&\hspace{-0.2cm}1 &\hspace{-0.2cm} 0 &\hspace{-0.2cm} 0&\hspace{-0.2cm} 0 &\hspace{-0.2cm} 0 \\
      \end{array}
    \right).}
$$

We now prove some properties of the obtained matrices.
\begin{lem}
{\rm The matrix $W_1$, $(W_2)$ obtained above, is a weighing matrix of weight four.}
\end{lem}
{\bf Proof.} We first prove that $W_1$ is orthogonal. Equivalently we prove that any two distinct rows of $W_1$ are orthogonal. We first treat rows of $R$ and $F$ separately. Let $F_i$, $F_j$ be the $i$th and $j$th rows of $F$ respectively ($i<j$). Then one of the followings hold.
\begin{itemize}
  \item The $1$ entries of $X_i$ and $X_j$ occur on different column. In this case, $F_i$, $F_j$ wouldn't have common columns of non-zero entries, which means they are orthogonal.
  \item $X_i(k)=X_j(k)=1$ for only one column number, say $1\leq k\leq m$. In this case, by the definition we have $F(i,2k-1)=F(i,2k)=1$, while $F(j,2k-1)=1$, $F(j,2k)=-1$. But for the other entries, $l\neq 2k-1,2k$ either of $F_i(l)$ or $F_j(l)$ are zero, hence we have $$F_i.F_j=F(i,2k-1)F(j,2k-1)+F(i,2k)F(j,2k)=0.$$
Thus $F_i$ and $F_j$ are orthogonal.
\end{itemize}

The orthogonality of row vectors of $R$ can be seen similarly. Now we prove that for any $i,j=1,\ldots,m/2$, the vectors $F_i$ and $R_j$, are orthogonal. The following two cases may occur.
\begin{itemize}
  \item $i=j;$ In this case, by the definition, we have $$R_i=[0,\ldots,0,F(i,2j_0-1),F(i,2j_0),0,\ldots,0,-F(i,2j_1-1),-F(i,2j_1),0,\ldots,0].$$
Where $j_0,j_1$ are the numbers of non-zero columns of $X_i$. Hence $F_i.R_i=F(i,2j_0-1)^2+F(i,2j_0)^2-F(i,2j_1-1)^2-F(i,2j_1)^2=1+1-1-1=0.$
  \item $i\neq j;$ In this case either $X_i$, $X_j$, and then $R_i$, $F_j$ have no common column having non-zero coordinates or $X_i$, $X_j$ share only one common column of non-zero elements. In the former case $R_i$, $F_j$ have zero inner product, indeed. In the case that $X_i$, $X_j$ share only one common column of non-zero elements, let $k$ be the index of the common column of non-zero entry. By the definition we will have $X_i(k)=X_j(k)=1$, so at first we have $F_j(2k-1)=F_i(2k-1),$ $F_j(2k)=-F_i(2k-1),$ now depending on the occurrence of $k$, one of the followings may occur:$$R_j(2k-1)=F_i(2k-1)\textrm{, }R_j(2k)=-F_i(2k-1), \textrm{ or }$$$$R_j(2k-1)=-F_i(2k-1)\textrm{, }R_j(2k)=F_i(2k-1).$$ In both cases the orthogonality of $F_i$ and $R_j$ follows.
\end{itemize}
Note that the matrix $X$ has two $1$ entries at each row and column, and in $W_1$, any entry $1$ of $X$ is replaced with $[1,1]$ or $[1,-1]$. Hence each row and column of $W_1$ has four entries equal to $\pm 1$, then $W_1$ has weight four, as desired. $\square$

Now we prove the semi-orthogonality of $W_1$, and $W_2$.
\begin{prop}
  {\rm The matrices $W_1$ and $W_2$ are semi-orthogonal. }
\end{prop}
{\bf Proof.} It suffices to prove for any appropriate $i,j$ the inner product of the $i$'th row of ${W_1}$ with the $j$'th row of ${W_2}$ is equal to $0,\pm 2$. Note that the matrices $X$ and $Y$ is chosen so that they never have equal rows. This implies that $X_i$ and $Y_j$ have at most one common $1$-entry. Hence the $i$'th row of ${W_1}$ and the $j$'th row of ${W_2}$ has either no common columns having non-zero coordinate or exactly two common columns having non-zero coordinates. In both cases depending on the corresponding entries, the inner product may be $0,$ or $\pm 2$. $\square$

The above statements result in the following theorem.
\begin{thm}
  {\rm Let $W_1$, and $W_2$ be two semi-orthogonal weighing matrices of order $m$. Then the signed graph with adjacency matrix $\mathcal{A}(W_1,W_2)$ has spectrum $[4^{m},-2^{2m}]$.}
\end{thm}
\begin{remark}
  {\rm Using Lemma \ref{con} and the mentioned construction of matrices $W_1$ and $W_2$ we may construct signed $8$-regular graphs on $3m$ nods, ($m$ even) having  spectrum $[4^m,-2^{2m}]$, for $m\geq6$. Using computer search we find the following instance with $12$ vertices.  

{\hspace{4.5cm}\begin{tikzpicture}[mystyle/.style={draw,shape=circle}]
\def\ngon{12}
\node[regular polygon,regular polygon sides=\ngon,minimum size=6cm] (p) {};
\foreach\x in {1,...,\ngon}{\node[mystyle] (p\x) at (p.corner \x){};}
\foreach\x in {1,...,\numexpr\ngon-1\relax}{
\foreach\y in {\x,...,\ngon}{
\draw (p1) -- (p5)[red];%
\draw (p1) -- (p11)[red];%
\draw (p1) -- (p8)[blue];%
\draw (p1) -- (p7)[red];%
\draw (p1) -- (p12)[red];%
\draw (p1) -- (p9)[red];%
\draw (p1) -- (p10)[red];%
\draw (p1) -- (p6)[red];%
\draw (p5) -- (p3)[blue];%
\draw (p5) -- (p11)[red];%
\draw (p5) -- (p12)[blue];%
\draw (p5) -- (p9)[blue];%
\draw (p5) -- (p10)[blue];%
\draw (p5) -- (p2)[red];%
\draw (p5) -- (p4)[red];%
\draw (p3) -- (p8)[blue];%
\draw (p3) -- (p11)[red];%
\draw (p3) -- (p6)[red];%
\draw (p3) -- (p7)[blue];%
\draw (p3) -- (p12)[red];%
\draw (p3) -- (p9)[blue];%
\draw (p3) -- (p10)[blue];%
\draw (p8) -- (p11)[red];%
\draw (p8) -- (p12)[red];%
\draw (p8) -- (p2)[red];%
\draw (p8) -- (p4)[blue];%
\draw (p8) -- (p9)[red];%
\draw (p8) -- (p10)[blue];%
\draw (p2) -- (p11)[blue];%
\draw (p11) -- (p4)[blue];%
\draw (p11) -- (p6)[blue];%
\draw (p11) -- (p7)[blue];%
\draw (p8) -- (p6)[blue];%
\draw (p6) -- (p2)[red];%
\draw (p6) -- (p4)[blue];%
\draw (p6) -- (p9)[red];%
\draw (p6) -- (p10)[blue];%
\draw (p7) -- (p12)[red];%
\draw (p7) -- (p4)[blue];%
\draw (p7) -- (p9)[blue];%
\draw (p7) -- (p10)[blue];%
\draw (p7) -- (p2)[blue];%
\draw (p12) -- (p4)[red];%
\draw (p12) -- (p2)[red];%
\draw (p9) -- (p2)[blue];%
\draw (p9) -- (p4)[red];%
\draw (p10) -- (p2)[red];%
\draw (p10) -- (p4)[blue];%
}}

\end{tikzpicture}}
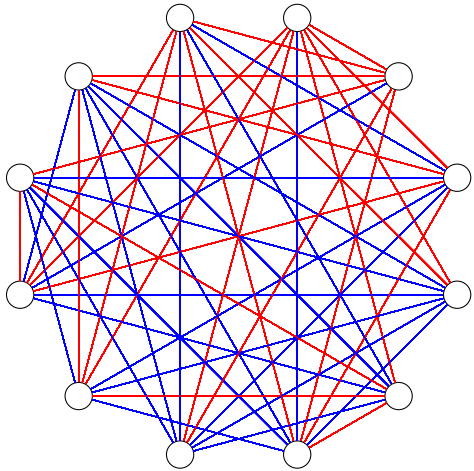
$$\textbf{Figure 2.}\textrm{ A signed graph with spectrum } [4^4,-2^8],$$
}
\end{remark}


\subsection{STE's with parameters $\mathbf{(0,\sqrt{k},-\sqrt{k})}$}
Signed graphs with spectrum $\sqrt{k}^{\frac{n}{2}},-\sqrt{k}^{\frac{n}{2}}$ are in correspondence with symmetric weighing matrices. Symmetric weighing matrices of small weights have been studied extensively and some constructions are announced, for more details see \cite{BKR,K,KS}. Recently Huang \cite{H}
has proved the Sensitivity Conjecture from theoretical computer science, by presenting a special signing of hypercubes. The obtained signed hypercubes admit spectrum $\pm\sqrt{k}$. His method is improved in \cite{BCKW} to a construction of symmetric weighing matrices. The following is their construction.

\begin{thm} \label{conswe} {\bf \cite{BCKW}} {\rm If $\Sigma$ is a signed $k$-regular graph of order $n$ with spectrum $[\sqrt{k}^{\frac{n}{2}},-\sqrt{k}^{\frac{n}{2}}]$, then the following is the adjacency matrix of a signed $(k+1)$-regular graph with spectrum $\sqrt{k+1}^{n},-\sqrt{k+1}^{n}$.$$Ac(\Sigma)=\left(
                                               \begin{array}{cc}
                                                 A_{\Sigma} & I_n \\
                                                 I_n & -A_{\Sigma} \\
                                               \end{array}
                                             \right).
$$}
\end{thm}
In \cite{DM} the authors have provided a family containing infinitely many signed $4$-regular graphs having just two distinct eigenvalues $2,-2$. More precisely for each $m\geq3$, they construct a signed $4$-regular graph with spectrum $[2^m,-2^m]$. We call their examples,\textit{ the $4$STE}'s. This yields the following result.
\begin{thm} \label{ram}
 {\rm For any $k\geq 4$ there are infinitely many connected signed $k$-regular graphs with distinct eigenvalues $\pm\sqrt{k}$.}
\end{thm}
{\bf Proof.} Note the graph $K_2$ has spectrum $\pm1$. Now repeatedly applying Theorem \ref{conswe}, for any positive $k$ we find a signed $k$-regular graph with distinct eigenvalues $\pm\sqrt{k}$. It is easily seen that the grounds of the signed graphs which are just constructed, are connected. As already discussed for any $m\geq3$, a signed $4$-regular graph with spectrum $[2^m,-2^m]$ exists. Let $\Sigma$ be a signed graph with spectrum $[2^m,-2^m]$. Then $Ac(\Sigma)$ would give a signed graph with spectrum $[\sqrt{5}^{2m},-\sqrt{5}^{2m}]$. By repeating the procedure we get a signed graph with spectrum $[\sqrt{6}^{4m},-\sqrt{6}^{4m}]$, $[\sqrt{7}^{8m},-\sqrt{7}^{8m}]$, $\ldots$, $[\sqrt{k}^{2^{k-4}m},-\sqrt{k}^{2^{k-4}m}]$. Since $m$ belongs to an infinite set of integer numbers hence the assertion follows. $\square$

The above theorem provides an infinite family of connected signed $k$-regular graphs having eigenvalues $\pm\sqrt{k}$. The ground of the mentioned graphs are bipartite. By the similar construction, starting with a symmetric conference matrix of order $6$, we provide other examples of signed $k$-regular graphs having eigenvalues $\pm\sqrt{k}$. In this case the ground is not bipartite. The following is the corresponding signed $5$-regular graph.

\hspace{5cm}\begin{tikzpicture}[mystyle/.style={draw,shape=circle}]
\def\ngon{6}
\node[regular polygon,regular polygon sides=\ngon,minimum size=4cm] (p) {};
\foreach\x in {1,...,\ngon}{\node[mystyle] (p\x) at (p.corner \x){};}
\foreach\x in {1,...,\numexpr\ngon-1\relax}{
\foreach\y in {\x,...,\ngon}{
\draw (p\x) -- (p\y)[blue];
\draw (p2) -- (p3)[red];
\draw (p2) -- (p4)[red];
\draw (p3) -- (p5)[red];
\draw (p4) -- (p6)[red];
\draw (p5) -- (p6)[red];}}
%
\end{tikzpicture}
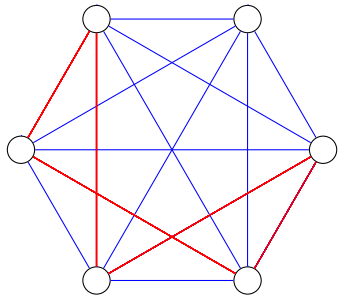
$$\textbf{ Figure 3. }\textrm{The Pentagon, its spectrum is }[\sqrt{5}^3,-\sqrt{5}^3] $$
Therefore we have the following result.
\begin{thm}
 {\rm For any $k\geq 5$ there is a signed $k$-regular graph on $3.2^{k-5}$ vertices with distinct eigenvalues $\pm\sqrt{k}$, in which the ground is not bipartite.}
\end{thm}
{\bf Proof.} By applying the method stated in Theorem \ref{conswe}, on the signed graph of Figure $3$, the assertion follows. $\square$

\subsection*{Concluding Remarks }
The vertex number of STE's which are provided at this paper are summarized at the following table. We proposed two infinite families of connected STE's. 

\hspace{4cm}\begin{tabular}{|c|l|}
                                                                                    \hline
                                                                                  $(t,\lambda_1,\lambda_2)$ & \textrm{ number of vertices of instances} \\\hline
                                                                                    $(1,3,-2)$ & $10$ \\
                                                                                   $(2,4,-2)$ & $15$, $6q$, $q\geq2$ \\
                                                                                    $(3,5,-2)$ & $21$ \\
                                                                                   $(0,\sqrt{k},-\sqrt{k})$ &$3.2^{k-4}$, $2^{k-3}m$, $m=3,4,\ldots$ \\

                                                                                    \hline
                                                                                  \end{tabular}

In \cite{BL} the authors have conjectured that every connected $k$-regular graph $\Gamma$ has a sign function $\sigma$ such that the maximum spectrum of the corresponding signed graph is at most $2\sqrt{k-1}$. Then they propose to construct new families of Ramanujan graphs by using lift of graphs. Our instance of signed $8$-regular graphs have maximum eigenvalue $4$, which is smaller than $2\sqrt{7}$. Hence an infinite family of $8$-regular Ramanujan graphs can be constructed. Also by Theorem \ref{ram} for any $k\geq 4$ an infinite family of signed graphs with maximum eigenvalue $\sqrt{k}$ is constructed, which again yields an infinite family of $k$-regular Ramanujan graphs.

Moreover it is well-known that the Kronecker product of two weighing matrices of orders $m$, $n$, and weights $\alpha,\beta,$ respectively, yields an orthogonal matrix of order $mn$, with weight $\alpha\beta$. Hence by the Kronecker product of the known weighing matrices with the above obtained matrices some new instances will be constructed.

\end{document}